\documentclass[11pt]{bmc_article}

\usepackage{cite} 					
\usepackage{url}  					
\usepackage{ifthen}  				
\usepackage{multicol}   		
\usepackage[utf8]{inputenc} 

\usepackage{algorithmic}
\usepackage{amsfonts,amsmath,amstext,epsfig,graphicx}
\usepackage{theorem}
\usepackage{verbatim,amssymb}
\usepackage[boxed]{algorithm2e}
\usepackage[tight,footnotesize]{subfigure}

\SetKwFor{Every}{every}{do}{end}

\newtheorem{remark}{Remark}

\newcommand{\midwor}[1]{\;\textrm{ #1 }\;}
\newcommand{\ds}{\displaystyle}

\renewcommand{\paragraph}[1]{\vspace{1ex}\textbf{#1}}

\renewcommand{\subset}{\subseteq}
\newcommand{\setB}{\mathcal{B}}

\newcommand{\setU}{\mathcal{U}}
\newcommand{\subsetU}{\mathcal{V}}

\newcommand{\ener}{\mathcal{E}}
\newcommand{\potu}{V}

\newcommand{\setC}{\mathcal{C}}

\newcommand{\pow}{P}

\newcommand{\att}{l}
\newcommand{\nois}{N}
\newcommand{\sinr}{\textrm{\textsc{sinr}}}

\newboolean{publ}

\newenvironment{bmcformat}{\baselineskip20pt\sloppy\setboolean{publ}{false}}{\baselineskip20pt\sloppy}

\begin{document}
\begin{bmcformat}


\title{Gibbsian Method for the Self-Optimization of Cellular Networks}

\author{Chung Shue Chen\correspondingauthor$^1$%
        \email{Chung Shue Chen\correspondingauthor - cs.chen@alcatel-lucent.com}
       and
        Fran\c{c}ois Baccelli$^2$%
        \email{Fran\c{c}ois Baccelli - francois.baccelli@ens.fr}}%
\address{%
    \iid(1)Network Technologies, Alcatel-Lucent Bell Labs, Centre de Villarceaux, 91620 Nozay, France\\
    \iid(2)Research group on Network Theory and Communications (TREC), INRIA-ENS, 75214 Paris, France}%
\maketitle

\begin{abstract}
In this work, we propose and analyze a class of distributed algorithms performing the joint optimization of radio resources in heterogeneous cellular networks made of a juxtaposition of macro and small cells.
Within this context, it is essential to use algorithms able to simultaneously solve the problems of channel selection, user association and power control. In such networks, the unpredictability of the cell and user patterns also requires distributed optimization schemes. The proposed method is inspired from statistical physics and based on the Gibbs sampler.
It does not require the concavity/convexity, monotonicity or duality properties common to classical optimization problems.
Besides, it supports discrete optimization which is especially useful to practical systems. We show that it can be implemented in a fully distributed way and nevertheless achieves system-wide optimality. We use simulation to compare this solution to today's default operational methods in terms of both throughput and energy consumption.
Finally, we address concrete issues for the implementation of this solution and analyze the overhead traffic required within the framework of 3GPP and femtocell standards.
\end{abstract}


\ifthenelse{\boolean{publ}}{\begin{multicols}{2}}{}

\section{Introduction\\} \label{sec:intro}

Today's cellular mobile radio systems strongly rely on highly
hierarchical network architectures that allow service providers to
control and share radio resources among base stations and clients in
a centralized manner.
With the foreseen exponentially
increasing number of users and traffic in the 4G and future
wireless networks, existing deployment and practice
becomes economically unsustainable.
Network self-organization and self-optimization
are among the key targets of future mobile
networks so as to relax the heavy demand of human efforts in the
network planning and optimization tasks and to reduce the system's
capital and operational expenditure (CAPEX/OPEX)
\cite{Socrates09,LTEbook2011,ICC10}. The next-generation mobile
networks (NGMN) are expected to provide a full coverage of broadband
wireless service and support fair and efficient radio resource
utilization with a high degree of operation autonomy and
intelligence.

Due to the emerging high demand of broadband service and new applications,
wireless networking also has to face the challenge of supporting fast increasing data traffic
with the requirement of spectrum and energy utilization efficiency \cite{Survey2011}.
To enhance the network capacity and support pervasive broadband service, reducing cell size is one of the most effective approaches.
Deployment of small cell base stations or femtocells has a great potential to improve the spatial reuse of radio resource and also enhance transmit power efficiency \cite{Femtocell2009}.
It is foreseen that the next generation of mobile cellular networks will consist of heterogeneous macro and small cells with different capabilities including transmit power and coverage range.
In such networks due to the unpredictability of the base station and user patterns, network self-organization and self-optimization become necessary.
Autonomic management and configuration of user association, i.e., assigning users to base stations, and radio resource allocation such as transmit power and channel selection would be highly desirable to practical systems \cite{3GPP2010R10PL}.

The primary objective of the present work is to design distributed algorithms performing radio resource allocation and network self-optimization for today's macro and small cell (e.g., 3GPP-LTE
\cite{LTEbook2011} and femtocell) mixed networks. In radio resource management, (i) power control, (ii) user association and (iii) channel selection are essential elements.
It is known that system-wide radio resource optimization is usually very challenging \cite{Luo2008}.
A joint optimization of user association, channel selection and power control is in general non-convex and difficult to solve, even if centralized algorithms are allowed \cite{ETT2011}.
Notice that in classical networks made of macro cells only, optimizing any of the above three elements independently can effectively improve the system performance. However, this may not be true in heterogeneous networks made of a juxtaposition of macro and small cells.  This would yield extra complexity and difficulties.
Besides, future wireless networks will typically be large, have fairly random topologies, and lack centralized control entities for allocating resources and explicitly coordinating transmissions with global coordination.
Instead, these networks will depend on individual nodes to operate autonomously and iteratively and to share radio resources efficiently.
We have to see how individual nodes can perform autonomously and support inter-cell interference management in a distributed way for finding globally optimal configurations.

To begin with, we give two examples to illustrate the problems that may happen when conducting these optimizations under macro and small cell networks, in both the downlink and uplink respectively.
Consider the downlink scenario in Figure~\ref{fig:eg1} where there are two mobile users $u$ and $v$ under the macro and small cell base stations (BS)
$a$ and $b$ which have different maximum transmit powers and coverage ranges.
Notice that user $u$ can be covered by the macro cell BS $a$ but it is located near the edge of $a$'s coverage. Meanwhile, it is too close to the small cell BS $b$ and this will have a strong impact on its received signal-to-interference-plus-noise-ratio (SINR). Here, transmit power optimization will not be effective without prior user association and channel selection optimization. One may consider the option in which users $u$ and $v$ both associate with the small cell $b$. However, this may overload BS $b$. From the viewpoint of load balancing, it is better to have the two users attached to different cells, e.g., user $u$ is attached to BS $a$. However, user $u$ will then have a low SINR as long as the two transmissions use a same channel. Clearly, one should consider assigning two different channels for these two transmitter-receiver pairs and hence conduct a joint user association and channel selection optimization with respect to the link characteristics of the possible combinations and their available channels. If  the system involves more users and cells, power control should be conducted as well to mitigate interference. This requires a joint optimization of all three elements.

Figure~\ref{fig:eg2} shows a similar problem in the uplink.
Consider that one first conducts user association optimization. Since user $v$ is closer to BS $b$ than to BS $a$, from the viewpoint of load balancing, the recommended user association should be as follows: user $u$ attaches to BS $a$ while user $v$ attaches to BS $b$. As user $u$ is far away from its BS $a$, the transmit power has to be high enough. This will however yield a strong interference to the signal received at BS $b$ which is transmitted from user $v$. Note that in this case, user association optimization, power control or even their joint optimization are not able to solve the problem. However,
if one also considers channel allocation and tries to select two different channels for these two transmitter-receiver pairs, a joint optimization will be able to resolve the conflict and enhance overall performance.


Let us now describe what aspects of the problem were considered so far
and the novelty of our approach. 
When each optimization is conducted separately,
the proper optimization sequence was studied
in \cite{SMARTA06,MDG07} for the 802.11 WLAN case,
based on careful experimental work and scenario analysis.
Explicit rules were proposed when the cell patterns have
a specific structure (e.g., in the hexagonal base station pattern case).
However, for situations where the cell and user patterns are unpredictable
as in the small cell case,
no simple and universal rule is known and a joint optimization
is necessary to achieve the best performance.



Various separate optimization problems were considered, mainly
under the assumptions of centralized coordination and
global information exchange.
For example the transmission powers
maximizing system throughput in the multiple interfering link case
leads to a non-convex optimization problem which
was studied in \cite{Chiang07a,ChenOien08}.
A power control algorithm that
guarantees strict throughput maximization in
the general SINR regime is reported in
\cite{MAPEL09}. It is built on multiplicative linear fractional programming,
which is used for optimization problems expressible as a difference of
two convex problems.
However, this algorithm requires a centralized control and is only efficient
for problem instances of small size due to the computation complexity.
There is a lack of efficient algorithm operating in a distributed manner and
ensuring global optimality in the above joint optimization.


Here, we propose and analyze a class of distributed algorithms performing
the joint optimization of radio resources in a generalized heterogeneous macro and small cell network.
Note that the optimization function does not have qualitative properties such as convexity or monotonicity.
The proposed solution is inspired from statistical physics and based on the Gibbs sampler
(see e.g., \cite{Geman84,bremaud99markov}).
It is a generalization of the work in \cite{ICC10}
which only takes into account power control and user association and is thus
limited to homogeneous mobile cellular networks.
The paper describes the algorithm,
shows that the latter can be implemented in a fully distributed
manner and nevertheless achieves
minimal system-wide potential delay,
reports on its performance, and
analyzes the overhead associated with the information exchange required in
the implementation of this solution in today's 3GPP-LTE and femtocell standards.
The rest of the paper is organized as follows.
Section~\ref{sec:model} describes the system model and problem setup.
Section~\ref{sec:Gibbs} presents the proposed solution.
Section \ref{sec:Simu} compares this solution to today's default operation in terms of throughput and energy consumption.
Section~\ref{sec:Overhead} investigates the overhead traffic generated by the algorithm.
Finally, Section~\ref{sec:conclusion} contains the conclusion.


\section{System Model and Problem Formulation\\} \label{sec:model}

We consider a reuse-1 cellular radio system with a set
$\setB$ of base stations serving a population
$\setU$ of users. For each user $u \in \setU$,
it is assumed that there is a pair of orthogonal
channels for the uplink and downlink.
We assume that there is no interference between the uplink and
downlink and we only consider the downlink.
However, the method can be generalized to the uplink as well.

We assume that users can associate with any
neighboring base station $b\in\setB$ in the network which could be a macro or small cell base station, which is
referred to as open access \cite{Femtocell2009}.
Today's default operation attaches each user $u$ to the base station
with the highest received power.
Note that this is clearly sub-optimal.
In general, if one simply associates users with the closest BS
or to that with the strongest received signal,
it is possible that some BSs have many users while others have only a few.
The resulting overload might lead to a degradation of the network capacity.

Let $\setC$ be the set of channels (e.g., frequency bands)
which are common to all base stations.
The base station serving user
$u$ is denoted by $b_u$ and is restricted to
some local set ${\mathcal B}_u$ of bases stations
(typically ${\mathcal B}_u$ is the set of
BSs the power of pilot signal of which is received
by user $u$ above some threshold).
The channel allocated by $b_u$ to user $u$ is denoted $c_u\in \setC$.
Here, for simplicity we consider that a user only takes one channel.
The transmission power used by base station $b_u$ to
$u$ is denoted by $\pow_u$.


The SINR at user $u$ is then:
\begin{eqnarray} \label{eq:sinr}
    \sinr_u = \frac
    {\pow_{u} \att(b_u,u,c_u)}
    {\nois_u(c_u) + \sum\limits_{v \in \setU, v \neq u}
    \alpha(b_u, b_v, c_u,c_v) \pow_v \att(b_v, u, c_v)}~,
\end{eqnarray}
where $\nois_u(c)$ denotes the thermal noise of user $u$ on channel $c$,
$\att(b_u,u,c)$ is the signal attenuation from BS $b_u$ to $u$ on channel $c$,
and $\alpha(b,b',c,c')$ represents the orthogonality factor
between
some user associated with BS $b$ on channel $c$
and some user associated with BS $b'$ on channel $c'$.

Note that it makes sense to assume that
$0 \leq \alpha(\cdot) \leq 1$ and that the
following symmetry holds: for all $b,b',c,c'$,
$$\alpha(b,b',c,c')=
\alpha(b',b,c',c).$$
Here are some examples: if adjacent channel interference
is negligible compared to co-channel interference, then one should take
$\alpha(b,b',c,c')=0$ for $c\ne c'$.
One may also assume that
$\alpha(b,b,c,c) = \alpha$ 
and $\alpha(b,b',c,c) = \beta$ for $b\ne b'$, 
where $\alpha$ and $\beta$ are some constants such that $\alpha < \beta$. 
The simplest case is that where $\alpha = \beta = 1$.


Under the additive white Gaussian noise (AWGN) model,
the achievable data rate at user $u$ in bit/s/Hz is given by:
\begin{eqnarray}
    \label{eq:rate}
    r_u = K \log(1+ \sinr_u)~,
\end{eqnarray}
where $K$ is a constant depending on the width of the frequency band.

To achieve network throughput enhancement while supporting bandwidth sharing fairness among users,
we adopt the notion of \textit{minimal potential delay fairness}
proposed in \cite{massoulie02bandwidth}.
This solution for bandwidth sharing is intermediate between max-min and proportional fairness.
It aims at minimizing the system-wide potential delay and is explained below.

Instead of maximizing the sum of throughputs, i.e., $\sum r_u$, which
often leads to very low throughput for some users,
we minimize the sum of the inverse of throughput, i.e., $\sum r_u^{-1}$, which can be seen as the
total delay spent to send an information unit to \textit{all} the users.
Note that minimizing $\sum r_u^{-1}$ penalizes very low throughputs.
More explicitly, a bandwidth allocation that provides
minimal potential delay fairness
is one that minimizes the following cost function:
\begin{eqnarray}
 \label{eq:PDF}
C =\sum_{u\in\setU} \frac{1}{r_u}~,
\end{eqnarray}
which is the network's aggregate transmission delay. It also indicates the long term
throughput that a user expects to receive from a fully saturated network.

For mathematical convenience (see below), in this paper, we minimize the cost function
\begin{eqnarray}
    \label{eq:energyf}
    \ener = \sum_{u\in\setU} \frac{1}{e^{\frac{r_u}K} -1}
= \sum_{u\in\setU} \frac{1}{\sinr_u}
\end{eqnarray}
instead of (\ref{eq:PDF}).
We call $\ener$ the global \emph{energy}, following the terminology of
Gibbs sampling.
Note that if one operates in a low SINR regime such
that the achievable data rate of a user is proportional
to its SINR, e.g., $r_u = K \sinr_u$, minimizing the potential delay $C$ is equivalent to minimizing the global energy $\ener$.

\begin{remark}
$\ener$ is a surrogate of $C$. We see that (\ref{eq:PDF}) and (\ref{eq:energyf})
have quite similar characteristics. The difference is that $(e^{\frac{r_u}{K}} - 1)^{-1}$ increases more significantly
than $r_u^{-1}$ when $r_u$ is low. As a result, the overall cost will increase more substantially.
So, minimizing $\ener$ rather than $C$ penalizes low throughputs more significantly
and favors a higher level of user fairness.
\end{remark}

By (\ref{eq:sinr}) and (\ref{eq:rate}), the global energy $\ener$ in (\ref{eq:energyf}) can be written as:
\begin{eqnarray}
    \label{eq:eneruCA}
\ener =
\sum_{u\in\setU} \frac{\nois_u(c_u) +
\sum\limits_{v \in \setU, v \neq u}
\alpha(b_u, b_v, c_u,c_v) \pow_v \att(b_v, u, c_v)}
{\pow_{u} \att(b_u,u,c_u)}
\end{eqnarray}
so that
\begin{eqnarray} \label{eq:enerUV}
    \ener = \sum_{u\in\setU} \frac{\nois_u(c_u)}
    {\pow_{u} \att(b_u,u,c_u)} + ~~~~~~~~~~~~~~~~~~~~~~~~~~~~~~~~~~~~~~~~~~~~~~~~~~~~~~~~~~ \nonumber \\
    \sum_{\{u,v\} \subset\setU}
    \left(  \frac{
    \alpha(b_u,b_v,c_u,c_v) \pow_v  \att(b_v, u, c_v)} {\pow_{u} \att(b_u,u,c_u)} + \right.
    \left.  \frac{
    \alpha(b_v, b_u, c_v, c_u) \pow_u  \att(b_u, v, c_u)} { \pow_{v} \att(b_v,v,c_v)} 
    \right).
\end{eqnarray}



The optimization problem consists in finding a configuration (also referred to as a \emph{state})
of user association, channel selection and power allocation which minimizes the above energy function.
It is clear that the problem has a high combinatorial complexity and is in general hard to solve for large networks.
However the additive structure of the energy can be used
to conduct its minimization using
a Gibbs sampler. This leverages the decomposition
of $\ener$ into a sum of local cost function for each user $u$ (say local energy $\ener_u$) which can be manipulated in a distributed way in the resource allocation. We explain this setup and optimization in the next section.


\section{Gibbs Sampler and Self Optimization\\}
\label{sec:Gibbs}

We now describe the distributed algorithm to perform the joint optimization of user association,
channel selection and power control.
It is based on a Gibbs sampler operating on a \emph{graph}
$\mathcal{G}$ of the network which can be defined as follows:
\begin{itemize}
  \item The set of \textit{nodes} in $\mathcal{G}$ is the set of users denoted by $u \in \setU$.
  \item Each node $u$ is endowed with a \emph{state} variable $s_u$ belonging to a finite set $\mathcal{S}$.
  The state of a node is a triple describing
 its user association, its channel and its
transmit power; this state denoted by 
  $s_u=\{b_u,c_u,P_u\}$.
  Here, we consider that transmit power is discretized.
  We denote the state of the graph by 
  $\mathbf{s} \triangleq (s_{u})_{u \in \setU}$.
  \item
Two user nodes $u$ and $v$ are \textit{neighbors} in this graph
if either (i) the power $P_0$ of the pilot signal received from a possible association base station for $v$ at $u$ is above some threshold, say $\theta$ or (ii) the power received from a possible base station for $u$
is above $\theta$ at $v$.
We denote the set of neighbors of $u$ by $\mathcal{N}_u$.
Notice that $v\in \mathcal{N}_u$ if and only if $u\in \mathcal{N}_v$.
\end{itemize}

Below, for all subsets $\subsetU\subset\setU$,
the cardinality of $\subsetU$ is denoted by $|\subsetU|$.

The global energy $\ener=\ener({\mathbf s})$ in 
(\ref{eq:enerUV}) derives from a
\emph{potential} function $\potu(\subsetU)$ \cite{bremaud99markov},
that is
\begin{equation}
    \label{eq:potential}
\ds \ener=
\sum_{\subsetU\subset\setU} \potu(\subsetU),
\end{equation}
where the sum bears on the set of all cliques of the graph defined
above and where the potential function $\potu(\cdot)$ has here
the following form:
\[
\left\{
\begin{array}{rll}
\potu(\subsetU) & = \ds
\frac{\nois_u (c_u)}{ \pow_{u} \att(b_u,u,c_u)}
& \hspace{-2.9cm} \midwor{if} \subsetU=\{u\},
\vspace{0.15cm}
\\
\potu(\subsetU) & = \ds
\frac{
\alpha(b_u,b_v,c_u,c_v) \pow_v  \att(b_v, u, c_v)}  { \pow_{u}  \att(b_u,u,c_u)} +
\frac{
\alpha(b_v,b_u,c_v,c_u) \pow_u  \att(b_u, v, c_u)}
 {\pow_{v}  \att(b_v,v,c_v)}
\midwor{~~if} \subsetU=\{u,v\},
\vspace{0.15cm}
\\
\vspace{0.15cm}
\potu(\subsetU) & = \ds
0
& \hspace{-2.9cm} \midwor{if}
|\subsetU| \geq 3.
\end{array}
\right.
\]



A global energy which derives from such a potential function satisfying the condition $\potu(\subsetU) = 0$ for $|\subsetU| \geq 3$ is hence amenable to a distributed optimization using the Gibbs sampler, which is based on
the evaluation of the local energy at each node:
\begin{eqnarray}
    \ener_u = \sum _{
\subsetU\subset\setU \ s.t.\  u\in\subsetU}
\potu(\subsetU) ~.
    \label{eq:local_energy}
\end{eqnarray}

Following the above definition of $\potu(\cdot)$, this can be re-written as:
\begin{eqnarray}
    \label{eq:loceneru}
\ener_u ({\mathbf s})
= \underbrace{
\frac{\nois_u (c_u)
+ \sum\limits_{v \neq u, v\in {\cal N}_u}
\alpha(b_u, b_v, c_u, c_v) \pow_v  \att(b_v, u, c_v)
}{ \pow_{u}  \att(b_u,u,c_u)}
}_{=1/(\sinr_u)} + \hspace{-0.2cm}
 \sum_{v\neq u, v\in {\cal N}_u} \hspace{-0.2cm}
 \frac{
\alpha(b_v, b_u, c_v, c_u) \pow_u  \att(b_u, v, c_u)}
 {\pow_{v} \att(b_v,v,c_v)}~. \hspace{-0.3cm}
\end{eqnarray}

The local energy can be written in the following form:
\begin{equation}
    \label{eq:funcP}
 \ener_u(\mathbf{s})= A_u(\mathbf{s}) + B_u(\mathbf{s}),
\end{equation}
where
$A_u(\mathbf{s})$ and $B_u(\mathbf{s})$ represent the first and second terms of (\ref{eq:loceneru}), respectively.
Notice that the first term $A_u(\mathbf{s})$ is equal to $1/\sinr_u$.
It is the ``selfish'' part of the energy function, which is small when
$\sinr_u$ is large. On the other hand, $B_u(\mathbf{s})$ is the ``altruistic'' part of the
energy, which is small when the power of the interference incurred
by all the other users 
because of $u$ is small compared to the power received from
their own base stations.

\begin{remark}
One can consider that $\ener_u$ consists of an individual cost of $u$ plus another term which corresponds to its impact on the others ($v \neq u$).
\end{remark}

\begin{remark}
The above formulation is meant to handle joint power, channel, and user association optimization.
However, it can easily be adapted to some special cases, e.g., to the case where the transmit power is a constant.
\end{remark}

In the following, we describe more precisely the Gibbs sampler
and its properties. 
First, we explain what it does.
Each BS separately triggers a state transition for one of
its users picked at random, say $u$,
using a local random timer.
This transition is selected based on the local energy $\ener_u$.
More precisely, given the state
$(s_v)_{v \neq u, \ v\in {\cal N}_u}$
of the neighbors of $u$, the new state $s_u$ is
selected in the set ${\mathcal S}_u$ of potential states
for user $u$
(this set is finite as power has been quantized to a finite set)
with the probability
\begin{eqnarray}
    \label{eq:distG}
\pi_u(s_u) = \frac{e^{- \frac{\ener_u(s_u, (s_v)_{v \in \mathcal{N}_u})}{T} }} {\sum\limits_{s \in \mathcal{S}_u} e^{- \frac{\ener_u(s, (s_v)_{v \in \mathcal{N}_u})}{T} }}~, ~~ s_u \in \mathcal{S}_u~,
\end{eqnarray}
where $T > 0$ is a parameter called the temperature.

We now list the properties of this sampler.
\begin{itemize}
\item These local random transitions drive the network to
a steady state which is the {\em Gibbs distribution} associated with
the global energy and temperature $T$, that is to a state
with the following distribution (in steady state):
\begin{eqnarray*}
 \pi_T(\mathbf{s})= c \cdot e^{-\ener(\mathbf{s})/T},
\end{eqnarray*}
with $c$ a normalizing constant. The proof is based on a reversibility
argument similar to that of \cite{bremaud99markov}.
\item This distribution puts more mass on low energy (small cost)
configurations and when $T \rightarrow 0$, the distribution $\pi_T(\cdot)$
converges to a Dirac mass at the state 
of minimal cost if it is unique (otherwise to a
uniform distribution on the minima).
\item This procedure is distributed in that the transition
of user $u$ only requires knowledge of the state of its neighbors.
We discuss the structure of message exchanges in more detail below.
\end{itemize}

The exact procedure
which users follow to conduct state transitions
is summarized in Algorithm \ref{alg:bstransition}.
Each user sets a timer, $t_u$, which decreases
linearly with time. We consider discrete time
in step of $\delta$ second(s) and simply set $\delta = 1$.
This timer has a duration randomly sampled
according to a geometric distribution.
When $t_u$ expires, a \emph{transition} of $u$ occurs
by which the state of this user
is updated as indicated above.


\begin{algorithm}[tbp]
\SetLine
\Every{$\delta$}{
\ForEach{$u$ 
}{
\eIf{$t_u \leq 0$}{
\ForAll{
$s$ in $\mathcal{S}_u$}{
$\ds 
\ener_u (s, (s_v, v\ne u))\leftarrow
A_u(s, (s_v, v\ne u)) + B_u(s,(s_v, v\ne u))
$\;
$d_u(s,(s_v, v\ne u)) \leftarrow\exp\left(
- \ener_u(s,(s_v, v\ne u)) / T \right)$\;
}
sample 
$s_u \in \mathcal{S}$ according to the probability
law $\pi_u(s, (s_v, v\ne u)) \triangleq d_u(s_u)/ \sum_{s \in \mathcal{S}_u} d_u(s, (s_v, v\ne u))$\;
sample $t_u\geq0$ with distribution $geom(1)$\;
}{$t_u\leftarrow t_u-\delta$\;}
}
}
\caption{State transition for the Gibbs sampler.}
\label{alg:bstransition}
\end{algorithm}

\subsection{A Few Remarks}

\paragraph{Greedy Variant}
One may consider to perform the state transition
by deterministically choosing the one that maximizes (\ref{eq:distG}) namely the \textit{best response} instead of selecting a state
according to the Gibbsian probability distribution.
It is known that a strategy of best response will drive the system to a
local minimum but not necessarily to an optimal solution.
Some discussions on the price of anarchy of a best response
algorithm can be found in \cite{Coucheney11} and references therein.  
The basic idea of the probabilistic approach described above 
is to keep a possibility to escape from being trapped in a local minimum.


\paragraph{Temperature and Speed of Convergence}
It is clear that the tuning of the temperature $T$ will
strongly impact the system's limiting distribution.
It has to be chosen by taking the tradeoff between the convergence speed and
the strict optimality of the limit distribution into account.

It is known that under conditions which ensure the compactness of the Markov forward operator and the irreducibility of the corresponding chain \cite{Liu95},
the Gibbs sampler will converge geometrically fast (for $T$ fixed) to
the Gibbs distribution.
In Section~\ref{sec:Simu}, 
we will present simulation results illustrating this convergence.

\paragraph{Annealed Variant}
For a fixed environment (i.e., user population, signal attenuation),
if one decreases $T$ as
$T  = 1/\ln(1 + t)$, where $t$ is time, then
the algorithm will drive the network to a state of minimal energy,
starting from any state.
A concrete proof of this result
is similar to that of \cite[pp.~311-313]{bremaud99markov}. This
proof is based on the notion of weak ergodicity of Markov chains
and reversibility argument and is omitted.

\subsection{Message Exchanges}

Two base stations, say $b$ and $b'$, 
are called \textit{implicit neighbors}
if there exist two neighboring users $u$ and $u'$
such that $u$ can associate to $b$ and $u'$ to $b'$, i.e., if
$b\in {\mathcal B}_u$,
$b'\in {\mathcal B}_{u'}$,
and either
$\alpha(b,b',c,c') P_0 \att(b',u,c') > \theta$
or $\alpha(b,b',c,c') P_0 \att(b,u',c) > \theta$ for some $c$, $c'$.
As we shall see, messages have to be exchanged between 
implicit neighbor base stations only (in addition to
those between users and their current association base station).

The necessity for message exchange comes from the need of sampling
$s_u$ in the algorithm. For this either user $u$ or its 
base station $\hat b_u$
before the sampling (below
we assume that the sampling takes place on $\hat b_u$)
has to have enough information to determine 
$\pi_u({s, (s_v, v\ne u)})$ or equivalently
$\ener_u(s, (s_v, v\ne u))$ for all $s\in {\mathcal S}_u$.
For this, some measurements and information exchange between
neighboring base stations and users are required.

The explicit definition of $\ener_u$ in (\ref{eq:loceneru}),
shows that for the evaluation of 
$A_u(s)$, a user $u$ will have to estimate the following data and
report them to its base station $\hat b_u$:
\begin{enumerate}
  \item the receiver noise: $N_u(c)$ on each channel $c$,
  \item the total received interferences:
  $\sum_{v\neq u} \alpha(b,b_v,c,c_v) P_v l(b_v,u,c)$, for each $c$ and for each $b \in {\mathcal B}_u$,
  and
  \item the path-loss or link gain: $l(b, u, c)$, for each $c$ and for each $b$ in the set ${\mathcal B}_u$.
\end{enumerate}

In order for $u$ or $\hat b_u$ to evaluate
$B_u(s)$, for all $s\in {\mathcal S}_u$,
each user $v\in {\cal N}_u$ will have to estimate
the following information and to report to its own base station $b_v$
(which will in turn communicate it to all its implicit neighbors
including $b_u$ on the backhaul network):
\begin{enumerate}
  \item the power of its received signal: 
  $P_v l(b_v,v,c_v)$, and
  \item the path-loss or link gain: $\alpha(b_v,b,c_v,c) l(b,v,c)$,
for each $c$ 
  and for each of $b\in {\mathcal B}_u$. 
\end{enumerate}

Note that the measurement of signal power, interference and path-loss
$l(b, u, c)$ for each considered channel from either its own base
station or neighboring base stations can be retrieved by the user
terminal from for example the measurement of available RSCP
(received signal code power) and/or RSSI (received signal strength indication). 

By the above information exchange, for each $u$, base station $\hat b_u$ is
able to compute $\ener_u$ for all $s\in {\mathcal S}_u$ and
hence to sample the new state $s_u$ of user $u$ according
to the above algorithm.
Notice that
inter-cell communication takes place
between implicit neighbor base stations only.
There is no need to transmit this
information via the wireless medium. We assume that this is
supported by the backhaul network. The amount of overhead
traffic generated by the algorithm can be evaluated. The results
on the matter are presented in Section~\ref{sec:Overhead}.



\section{Simulation and Comparison\\}
\label{sec:Simu}

A performance investigation of the proposed solution is conducted below.
We implement Algorithm~\ref{alg:bstransition} and compare its performance with today's 3GPP default operations \cite{3GPPTR2010} by discrete event simulations.

In the current standard and 3G implementations,
base stations are usually configured with a nominal fixed transmission power such that the pilot signal can be received by terminals over the covered area. The downlink transmit power is often the maximum allowable power as well for a better user reception and coverage. Note that the pilot signal is broadcasted continuously to allow user equipments (UE) to perform channel measurements and appropriate tuning.
In user association, the current practice consists in attaching a user to the BS received with the strongest signal strength (rather than the nearest base station). Note that this could lead to attaching the users to a far macro cell BS which has a higher transmit power than that of a nearer small cell BS. This is in general sub-optimal.
In channel allocation, the current practice often follows a heuristic scheme where channels of a BS are assigned to its users simply in a round-robin fashion, i.e., sequentially, and in such a way that the numbers of users on each channels are well balanced and almost equal.



In the simulations, we consider that mobile users are uniformly distributed in a geographic area of $1000$ meters times $650$ meters and we adopt the 3GPP-3GPP2 spatial channel model \cite{3GPPChannelModel07}. The distance dependent path-loss is given by:
\begin{eqnarray}
    l^{\textrm{(dB)}}(d) = - 30.18 - 26 \log_{10}(d) - X_{\sigma}^{\textrm{(dB)}}~,
\label{eqn:PLeqn}
\end{eqnarray}
where $d$ is the transmitter-receiver distance and $X_{\sigma}$ refers to log-normal shadowing with zero mean and standard deviation 4 dB.
With operating temperature 290 Kelvin and bandwidth 1 MHz,
the thermal noise $N_u$ is equal to $4.0039 \times 10^{-15}$ W, for all $u$.

Here we consider that there are two macro cell base stations with fixed locations as shown
in Figure~\ref{fig:BS} and a number of small cell base stations which are randomly located
in the geographical area.
The maximum transmit power 
of macro and small cell base stations are 40W and 1W respectively.
We assume that $\pow_\delta = 0.1$ W. 
In the simulation,
we consider a simple system where $\alpha = 1$ and each user only takes one channel.

\subsection{Numerical Examples\\}

To begin with, we illustrate the effectiveness of the algorithm by some examples with randomly generated small cell BS and users, as shown in Figures~\ref{fig:beforeoptimiz}--\ref{fig:afteroptimiz}.
To have readable graphical representation and comparison
of the user association, channel allocation and transmission power before and after optimization,
in these examples, we consider that the path-loss is simply distance dependent without log-normal shadowing.
So, a user who is farther from a BS has a larger path-loss due to the larger distance.
A line connecting a BS and a user indicates the user association and its thickness represents the
strength of the transmit power. In these examples, we consider that there are two orthogonal channels in each BS,
which are represented by different colors and line styles.

Our simulations 
show that the proposed solution significantly outperforms
the by-default configuration in both system throughput (in b/s/Hz) and power consumption efficiency (in b/s/Hz/W).
Note that the latter has been improved by several orders of magnitude (also because our representation
of the default operation has no power control mechanism).
Figure~\ref{fig:num_cong} shows the corresponding convergence of the algorithm in the above three examples.
We see that the algorithm usually converges in a few hundreds of iterations and is hence practical.

\subsection{Average Performance\\}


Secondly, we compare the performance of the proposed optimization with the default operation,
with a fixed number of 32 BS (including the two macro BS) but with different
numbers of users (denoted by $M$), i.e., different user densities, and different numbers of orthogonal
channels (denoted by $K$). Users and small cells are randomly generated in the geographical
area. For each $(M, K)$, 500 different topologies are sampled and the
performance metrics are then averaged out.

Table~\ref{tab:throughput} shows the
the enhancement of the system throughput and of the power efficiency
obtained by the joint optimization. 
Observe that for a given $M/K$ ratio, the spectrum utilization efficiency
that results from the optimization increases with $K$.
This observation is important for e.g., in 3GPP HSDPA
(High Speed Downlink Packet Access) and LTE, where a high number
of users and a high number of resources are typical.

\section{Evaluation of Overhead Traffic\\}
\label{sec:Overhead}

The aim of this Section is to evaluate the overhead traffic
generated by the algorithms in a specific scenario which
is based on the assumption that nodes form realizations
of Poisson point processes in the Euclidean plane. These
assumptions allow us to use elementary stochastic geometry
to get estimates of this overhead traffic.

We concentrate on the channel selection and power control optimization,
when assuming that users are associated with their closest or best base station.
The overhead traffic has two main components: (i) the uplink radio traffic
and (ii) the backhaul traffic.

\subsection{Setting\\}

The uplink radio overhead traffic
is comprised of the set of messages that are sent by each mobile
to its serving base station and that inform the latter of the path-loss
that it experiences from each of its neighboring base stations. These data
are required to run the algorithm, see e.g., (\ref{eq:loceneru}).
If one denotes by $\tau$ the frequency of the beaconing signals
from the base stations and if one assumes that the users report
their path-loss variables at each beacon, each mobile has to
report $N \times \tau$ path-loss per second when the number of its neighboring base stations
is $N$.

On the other hand,
the backhaul traffic is between base stations (it is typically transported
by a wireline infrastructure). We will say here that two base stations
are neighbors if one of them has customers which see the other as a neighboring base station.

Consider a pair of neighboring base stations.
Let $M_1$ denote the number customers of the first base station (say BS 1) which see the second
(say BS 2) as a neighboring base station. Let $M_2$ be the symmetrical variable. Then
the global backhaul traffic between the two stations is given by:
\begin{equation}
    \sum_{i=1}^{M_1} \tau N_{1, i}
+ \sum_{j=1}^{M_2} \tau N_{2, j}
\end{equation}
where $N_{1, i}$ denotes the number of neighboring base stations of BS 1 for user $i$ and
$N_{2, j}$ denotes the number of neighboring base stations of BS 2 for user $j$.
Note that their definitions are symmetric.


\subsection{Stochastic Geometry Model\\}

We first describe the model for the overhead traffic for a purely macro cellular
network and then for an heterogeneous network with both macro and small cells.

\subsubsection{Macro Cell Model}
The base stations are assumed to form a Poisson point process of
intensity $\lambda_m$ in the Euclidean plane.
The users are assumed to form an independent Poisson point process of
intensity $\lambda_u$ in the Euclidean plane.
The association of the users to the closest BS
makes the association region of a base station to be the Voronoi
cell of this base station with respect to the collection of base stations.
This association together with the downlinks are depicted in Figure~\ref{fig1}.

The mean number of users of a typical cell, 
denoted by $\overline M$, is equal to $\lambda_u / \lambda_m$. 
In our model, we will assume that all users in a cell have for neighboring
base stations the Delaunay neighbors of the base station which is the nucleus
of the cell. This is depicted in Figure~\ref{fig2}.

The mean number of Delaunay neighbors of a typical node is 6
and its coefficient of variation
$CV(N) = \sqrt{Var(N)}/E(N)$ is $CV(N) = 0.222$ (see e.g., \cite{Moller}).

Hence, a rough estimate of the mean
uplink radio overhead traffic is:
\begin{equation}\label{eq:uplink_overhead}
    \overline R \approx 6\tau \frac{\lambda_u}{\lambda_m}.
\end{equation}
This is only an estimate because there is a correlation between the number
of users in a cell and the number of neighbors of the nucleus of this cell.
We now give an upper bound on $\overline R$ in complement of this estimate.

The second moment of the number of users in a cell is (see \cite{Moller}):
\begin{equation}
    E(M^2)= \frac{\lambda_u}{\lambda_m} +1.280 \frac{\lambda_u^2}{\lambda_m^2}.
\end{equation}

The second moment of the number of neighbors of a cell is given by:
\begin{equation} 
    E(N^2)= Var(N) + E(N)^2 = 37.7742.
\end{equation}

One can then use the Cauchy-Schwarz inequality to get the following upper-bound:
\begin{equation}
\overline R \leq\tau
\sqrt{\left(\frac{\lambda_u}{\lambda_m} +1.280 \frac{\lambda_u^2}{\lambda_m^2}\right)
E(N^2)}~.
\end{equation}

%
%

Consider now a typical backhaul link, namely a typical Delaunay edge.
A rough estimate of the mean backhaul overhead traffic on this link is then given by:
\begin{equation}
    \overline B = 2 \overline R \approx 12 \tau \frac{\lambda_u}{\lambda_m}.
\end{equation}
The Cauchy-Schwarz inequality can again be used to get an upper-bound.

\subsubsection{Macro and Small Cell Model}

In this section, we assume that each small
cell has a radius of coverage and that all users covered by the
small cell are attached to it. We also assume that
small cell rarely overlap. The users not covered by a small cell
are attached to the closest macro base station. This is depicted in
Figure~\ref{fig3}.

We assume that the small cell base stations form an independent point process
of intensity $\lambda_s$ and that the radius of coverage is $\rho$.
The mean number of users in a small cell is thus given by:
\begin{equation}
    \overline M_s = \lambda_u \pi \rho^2
\end{equation}
while the mean number of users attached to a macro cell is given by:
\begin{equation}
    \overline M_m=  \frac{\lambda_u}{\lambda_m} -\lambda_u \lambda_s \pi \rho^2.
\end{equation}

This formula is only valid under 
that the Boolean model with intensity $\lambda_s$ and
radius $\rho$ has only rare intersections of balls.

We declare neighbors of a macro cell its macro cell neighbors, defined as above,
and all small cells whose base station is located in the macro
cell in question or in one of its neighboring macro cells.

We declare neighbors of a small cell the base station of the
macro cell it is located in and the macro neighbors of the latter
as well as the small cells located in these macro cells.

Since the mean number of small cells per macro cell
is $\frac{\lambda_s}{\lambda_m}$, the mean number of small cells
neighbors of a macro cell is:
\begin{equation}
    \overline N_m^s = 7\frac{\lambda_s}{\lambda_m},
\end{equation}
while the mean number of macro cells neighbor of a macro cell is still 6.

The mean number of macro cells neighbors of a small cell is 7
and the mean number of small cells neighbor of a small cell is:
\begin{equation}
    \overline N_s^m = 7\frac{\lambda_s}{\lambda_m}~.
\end{equation}

Thus,
the mean uplink radio overhead traffic on a macro cell is given by:
\begin{eqnarray}
\overline R_m  & \approx & 6 \tau \overline M_m + \overline N_m^s \overline M_s \nonumber \\
& \approx &
6 \tau \left(\frac{\lambda_u}{\lambda_m}-\lambda_u \lambda_s \pi \rho^2\right)
+7\tau \frac{\lambda_s}{\lambda_m}\left(\lambda_u \pi \rho^2\right)
\end{eqnarray}
whereas that on a small cell is given by:
\begin{eqnarray}
\overline R_s  & \approx & 7 \tau \overline M_m + \overline N_s^m \overline M_s \nonumber \\
& \approx &
7 \tau \left(\frac{\lambda_u}{\lambda_m}-\lambda_u \lambda_s \pi \rho^2\right)
+7\tau \frac{\lambda_s}{\lambda_m}\left(\lambda_u \pi \rho^2\right) .
\end{eqnarray}

The  mean backhaul traffic on a link between two macro base stations
is $2 \overline R_m$, whereas that between a macro base station and a small base
station is equal to $\overline R_m+\overline R_s$.

These mean values can be complemented by bounds
using second moments.


\section{Conclusion\\}
\label{sec:conclusion}

In this paper, we analyzed the problem of radio resource allocation in heterogeneous cellular networks composed of macro and small cells with unpredictable cell and user patterns.
To solve the problem, we proposed a joint optimization of channel selection, user association and power control.
The proposed solution, which is based on the Gibbs sampler,
is implementable in a distributed manner and nevertheless
achieves minimal system-wide potential delay,
regardless of the initial state.
We investigated its performance and estimated the expected overhead.
Simulation result and comparison to today's default operations
have shown its high effectiveness in terms of energy consumption.
Because of its operational simplicity,
this distributed optimization
approach is expected to play an important role
in the future of heterogeneous wireless networks.



\section*{Acknowledgments}
  \ifthenelse{\boolean{publ}} {\small}{}
  The work presented in this paper has been carried out at LINCS (www.lincs.fr) and 
  under the INRIA-Alcatel-Lucent Bell Labs Joint Research Center.
  A part of this work was presented in \cite{VTC11} at the IEEE VTC workshop on Self-Organizing Networks.
	We would like to thank Laurent Thomas, Laurent Roullet and Vinod Kumar of Alcatel-Lucent Bell Labs for their valuable discussion and continuous support to this work.



{\ifthenelse{\boolean{publ}}{\footnotesize}{\small}
 \bibliographystyle{bmc_article}  
  \bibliography{bmc_article} }     

\ifthenelse{\boolean{publ}}{\end{multicols}}{}




\newpage


\section*{Table:}

\begin{table}[ht]
    \centering
\renewcommand{\arraystretch}{1.2}
    \caption{User average throughput: b/s/Hz,~ Power efficiency: b/s/Hz/W}
     \small{
\begin{tabular}
    {c|c|c|c}
    \hline
 & Default Operation & After Optimization & Performance Gain (times)\\  
\hline \hline
\hspace{-0.1cm}$M$ = 32,  $K$ = 1 & 0.245, 0.0143  & 1.216, 1.937  & 4.96, 135 \\ \hline
\hspace{-0.1cm}$M$ = 64,  $K$ = 2 & 0.312, 0.0186  & 1.583, 2.685  & 5.07, 144 \\ \hline
\hspace{-0.1cm}$M$ = 96,  $K$ = 3 & 0.356, 0.0210  & 1.829, 3.149  & 5.14, 150 \\ \hline
\hspace{-0.1cm}$M$ = 160, $K$ = 5 & 0.368, 0.0228  & 1.973, 3.488 & 5.36, 153 \\ \hline
\end{tabular}}
    \label{tab:throughput}
\end{table}


\section*{Figures:}

\begin{figure}[ht]
\begin{center}
    \includegraphics[width=0.55 \linewidth]{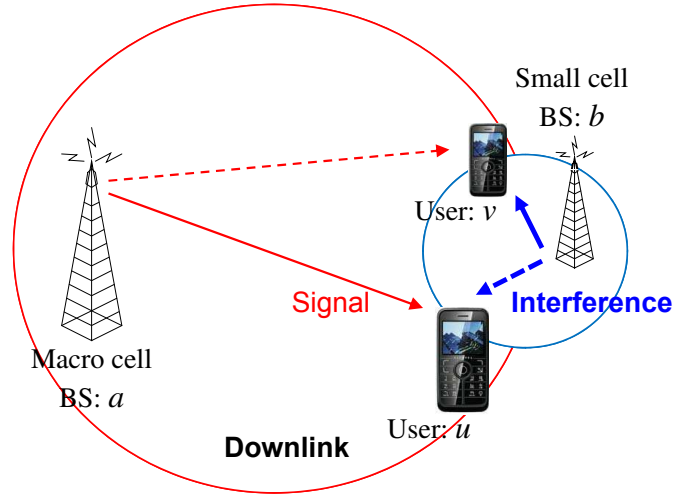}
\caption{Since user $u$ is far from its BS $a$, the received signal at user $u$ may suffer strong interference due to the transmission of small cell BS $b$ destined to user $v$.}
    \vspace{0.1cm}
\label{fig:eg1}
\end{center}
\end{figure}

\begin{figure}[ht]
\begin{center}
    \includegraphics[width=0.55 \linewidth]{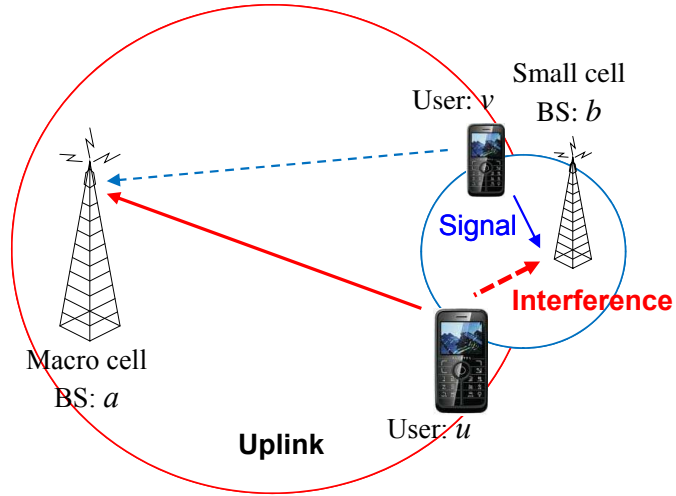}
        \vspace{0.1cm}
\caption{The signal received at BS $b$ sent from user $v$ can be strongly interfered by the transmission of user $u$ since $u$ has to use a relatively high power in order to send its signal to BS $a$ in long distance.}
\label{fig:eg2}
\end{center}
\end{figure}

 \begin{figure}[ht]
\centering
\includegraphics[width=0.65 \textwidth]{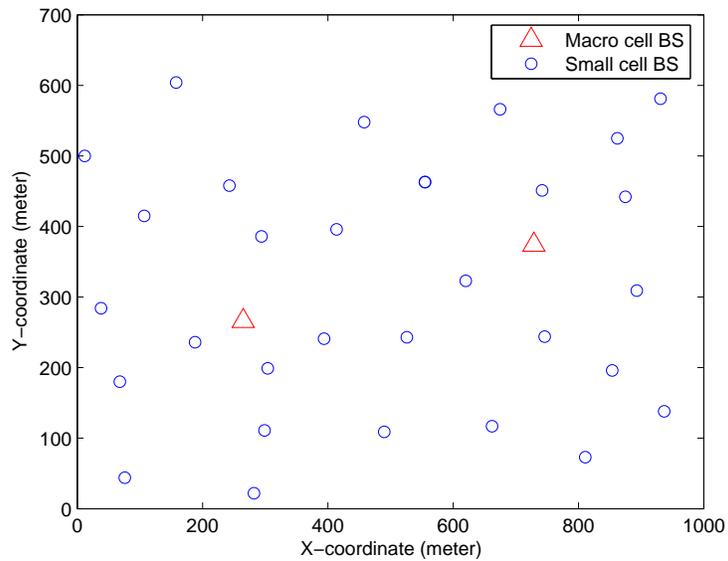}
\caption{The geographic location of macro and small cell base stations (example).}
\label{fig:BS}
\end{figure}


    \begin{figure*}[htb]
    \centering
\hspace{-0.1cm}
        \subfigure[Example~1: i) 8.7 b/s/Hz, ii) 0.012 b/s/Hz/W  ~~~~~~~~~~~~~~~~~~~~~~~~~~~~~~~~~]
        {
            \includegraphics[width=0.45\textwidth]{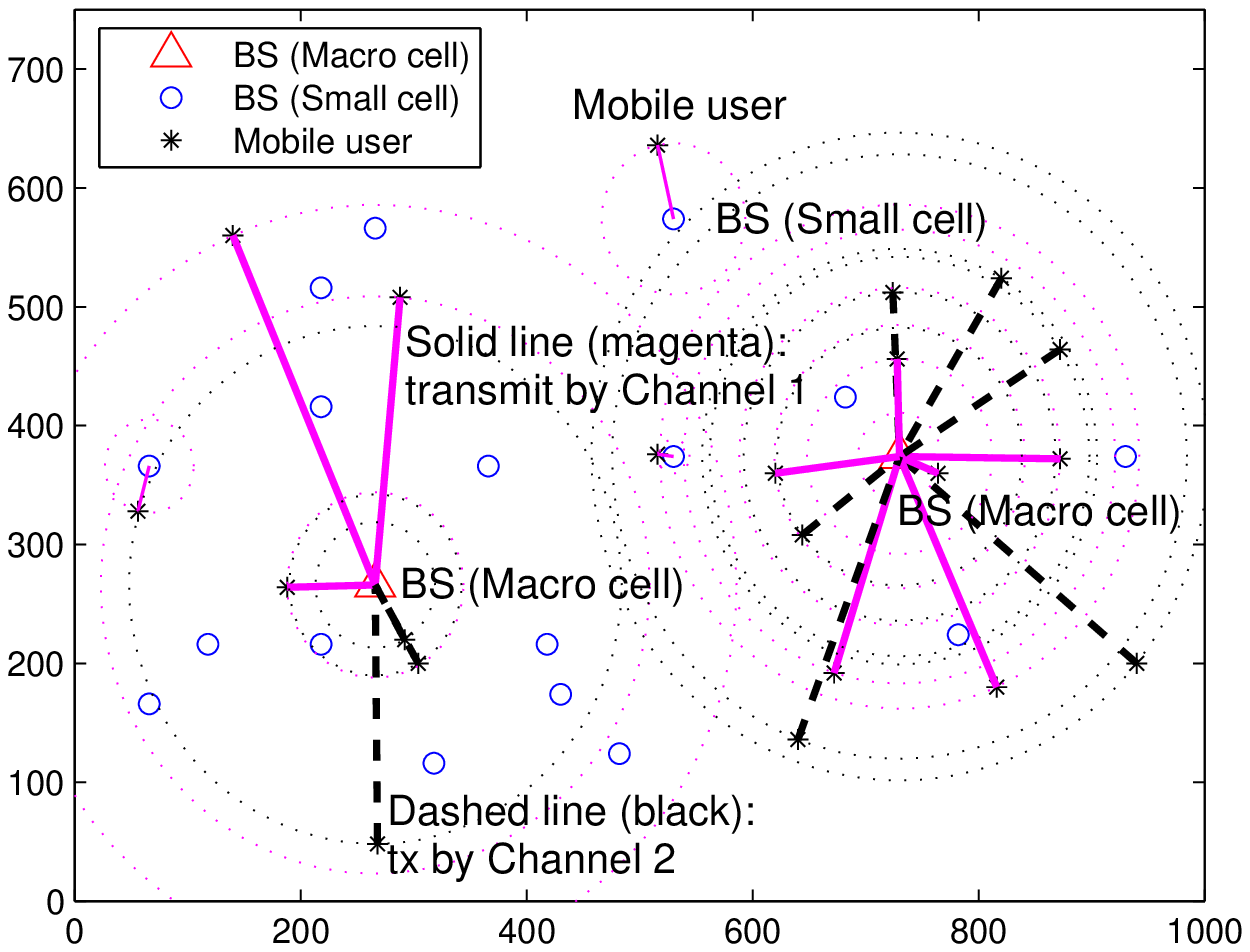}
        }
\hspace{0.3cm}
        \subfigure[Example 2: i) 35 b/s/Hz, ii) 0.106 b/s/Hz/W]
        {
            \includegraphics[width=0.45\textwidth]{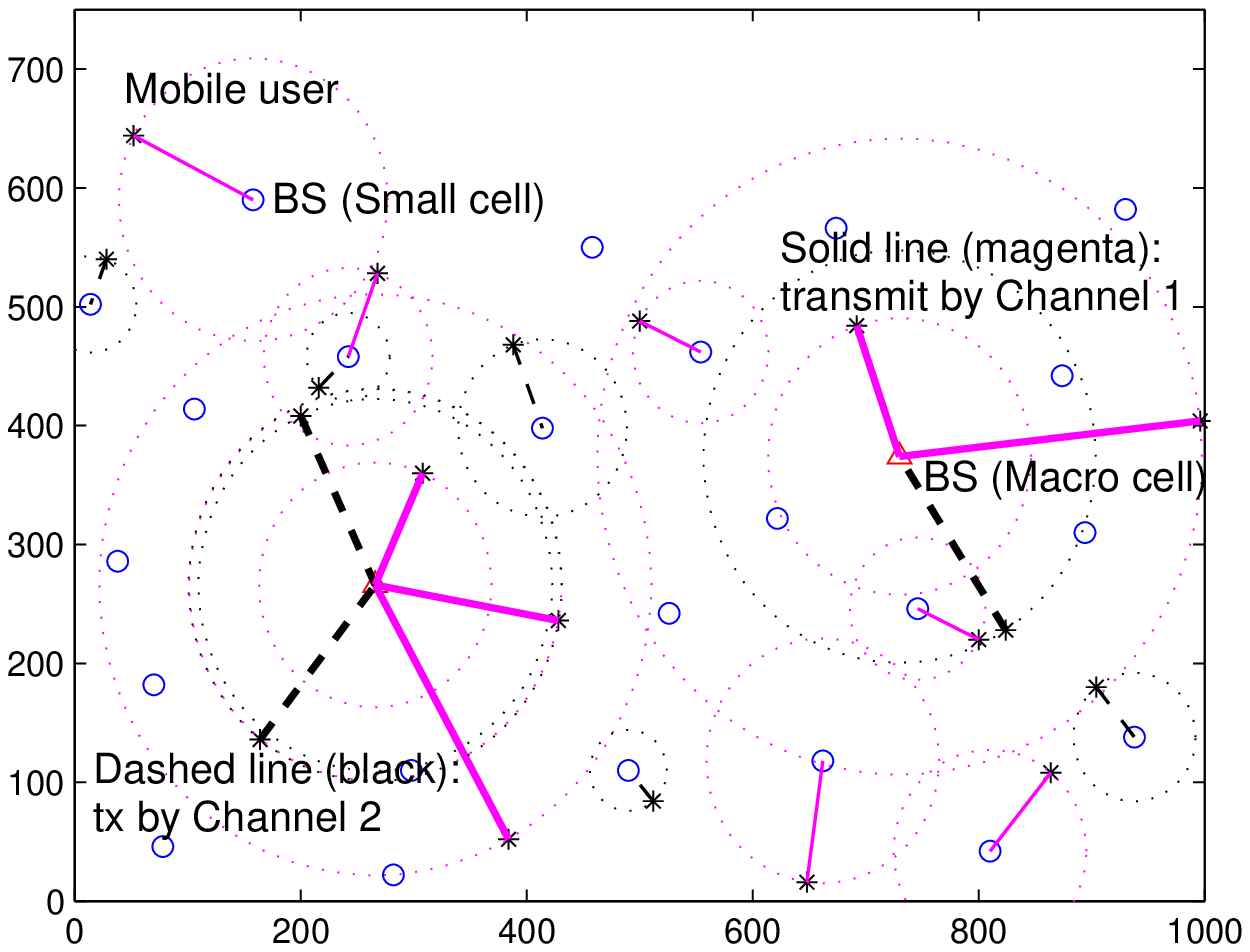}
        }
\hspace{0.0cm}
        \subfigure[Example 3: i) 7.5 b/s/Hz, ii) 0.009 b/s/Hz/W]
        {
        		\includegraphics[width=0.45\textwidth]{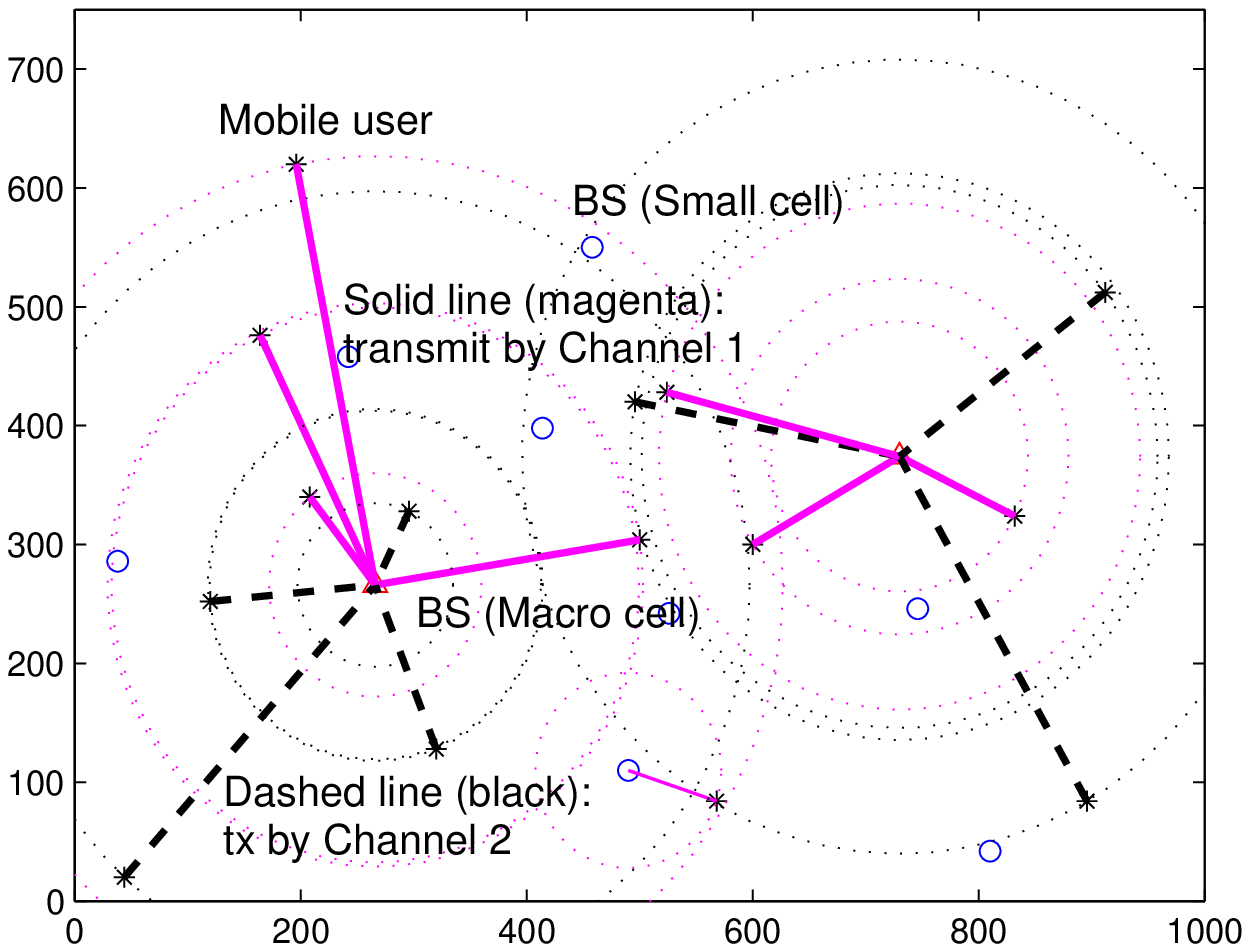}
        }
        \caption{Network before optimization (default operation). (a) Example 1: users are concentrated and fewer than BS. (b) Example 2: users are distributed and fewer than BS. (c) Example 3: more users than BS. Performance measure: i) system throughput, and ii) power efficiency. There are two orthogonal channels represented by solid-magenta and dashed-black lines.}
        \label{fig:beforeoptimiz}
    \end{figure*}


    \begin{figure*}[htb]
    \centering
\hspace{-0.1cm}
        \subfigure[Example 1: i) 43.5 b/s/Hz, ii) 3.45 b/s/Hz/W ~~~~~~~~~~~~~~~~~~~~~~~~~~~~~~~~~]
        {
            \includegraphics[width=0.45\textwidth]{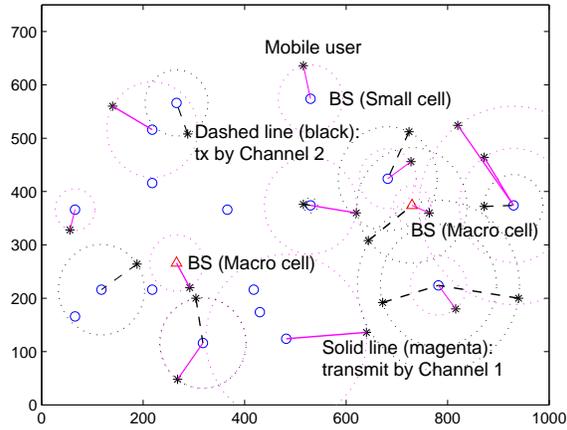}
        }
\hspace{0.3cm}
        \subfigure[Example 2: i) 75 b/s/Hz, ii) 4.41 b/s/Hz/W]
        {
            \includegraphics[width=0.45\textwidth]{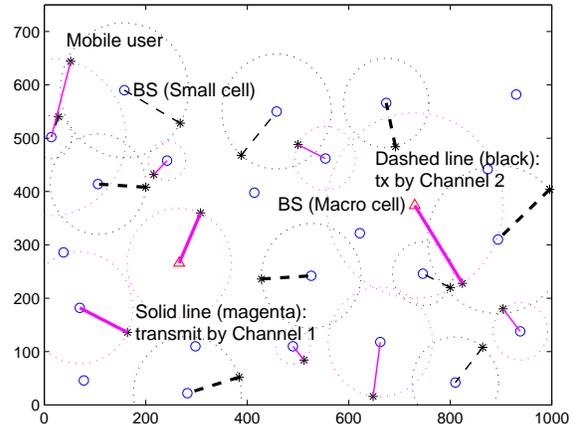}
        }
\hspace{0.0cm}
        \subfigure[Example 3: i) 32 b/s/Hz, ii) 2.13 b/s/Hz/W]
        {
            \includegraphics[width=0.45\textwidth]{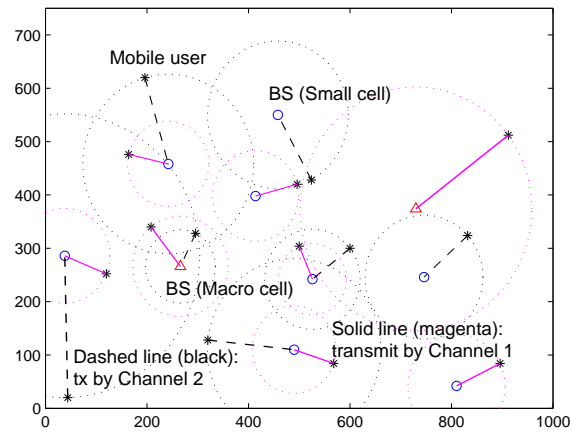}
        }
        \caption{Network after proposed joint optimization. Both the system throughput (b/s/Hz) and power utilization efficiency (b/s/Hz/W) are significantly improved.}
        \label{fig:afteroptimiz}
    \end{figure*}

\begin{figure}[htbp]
\centering
\includegraphics[width=0.65 \textwidth]{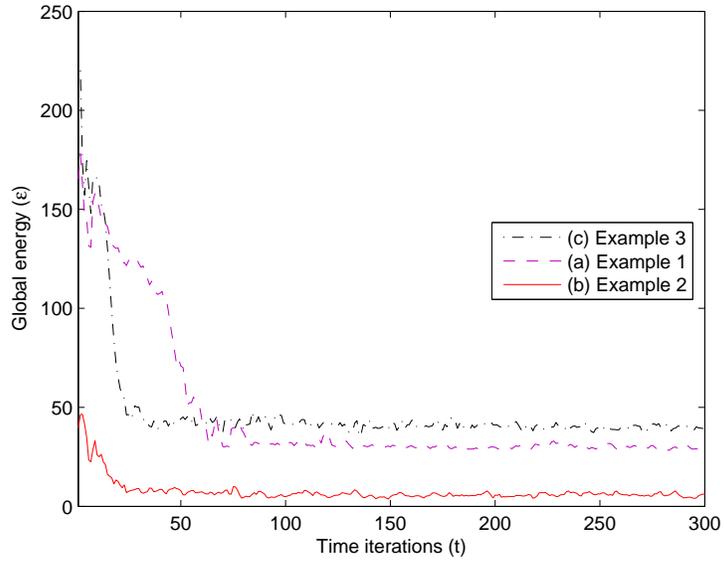}
\caption{Convergence of the algorithm: (a) Example 1, (b) Example 2, and (c) Example 3, respectively.}
\label{fig:num_cong}
\end{figure}

\begin{figure}[htb]
    \vspace{0.3cm}
\centering     \hspace{0.5cm}
\includegraphics[width=0.5\textwidth]{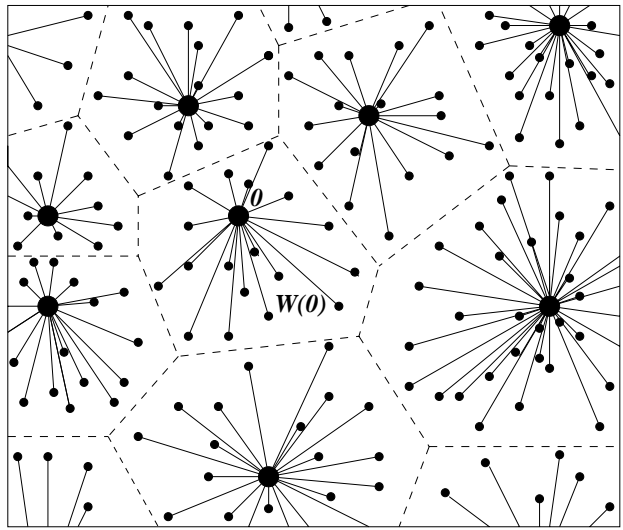}
    \vspace{0.2cm}
\caption{The dashed lines represent the boundaries of the cells. The solid lines
link from the base stations to the users which they serve.}
    \vspace{0.4cm}
\label{fig1}
\end{figure}

\begin{figure}[htb]
    \vspace{0.3cm}
\centering
\includegraphics[width=0.5\textwidth]{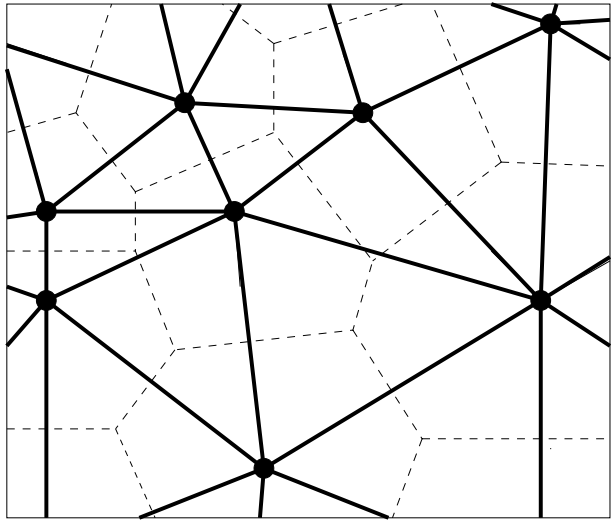}
    \vspace{0.2cm}
\caption{The solid lines
represent the Delaunay graph and serve as model for the backhaul network.}
\label{fig2}
\end{figure}

\begin{figure}[htb]
    \vspace{0.3cm}
\centering
\includegraphics[width=0.5\textwidth]{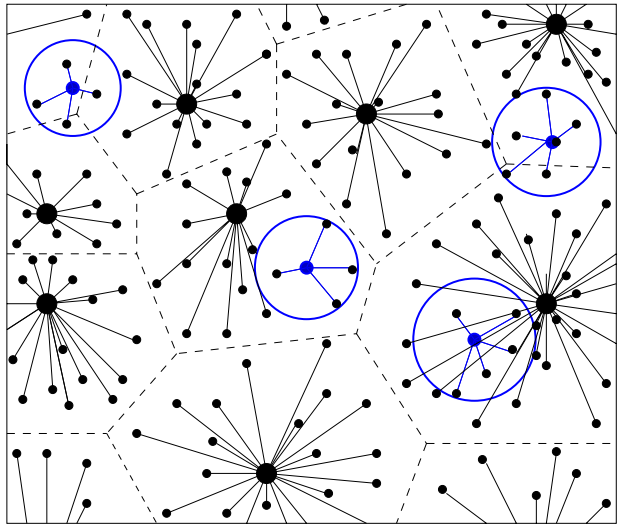}
    \vspace{0.2cm}
\caption{The discs represent the small cells. The solid lines again represent
the links from the base stations to the users they serve.}
\label{fig3}
\end{figure}

\end{bmcformat}
\end{document}